\newcommand{\R}{I\!\!R}
\newcommand{\Z}{\mathbb Z}
\newcommand{\N}{I\!\!N}
\newcommand{\Ker}{\mathrm{Ker}}
\newcommand{\Img}{\mathrm{Im}}
\newcommand{\Spl}{\mathrm{Sp}}
\newcommand{\spl}{\mathrm{sp}}
\newcommand{\GL}{\mathrm{GL}}
\newcommand{\Lin}{\mathrm{Lin}}
\newcommand{\Bil}{\mathrm{Bil}}
\newcommand{\Bils}{\mathrm{Bil}_{\mathrm{sym}}}
\newcommand{\Dim}{\mathrm{dim}}
\newcommand{\Id}{\mathrm{Id}}
\newcommand{\mul}{\mathrm{mul}}
\newcommand{\Dcal}{{\mathcal D}}
\newcommand{\Hcal}{{\mathcal H}}
\newcommand{\iMaslov}{\mathrm i_{\mathrm{Maslov}}}
\newcommand{\iMaslovred}{\mathrm i_{\mathrm{Maslov}}^{\scriptscriptstyle\mathrm{red}}}

\newcommand{\Pcal}{\mathcal P}
\newcommand{\Scal}{\mathcal S}
\newcommand{\Kcal}{\mathcal K}
\newcommand{\dgn}{\mathrm{dgn}}
\newcommand{\sgn}{\mathrm{sgn}}
\newcommand{\dd}{\mathrm d}
\newcommand{\red}{{\mathrm{red}}}
\newcommand{\mathbfB}{B^{\int}}
\newcommand{\gfrak}{\mathfrak g}
\newcommand{\Rcal}{\mathcal R}
\newcommand{\Ncal}{\mathcal N}
\newcommand{\Ycal}{\mathcal Y}
\newcommand{\Hess}{\mathrm{Hess}}
\newcommand{\vfrak}{\mathfrak v}
\newcommand{\wfrak}{\mathfrak w}
\newcommand{\II}{\mathrm I\!\mathrm I}
\newcommand{\cte}{\mathrm{Const}}

\documentclass[oneside,draft,letterpaper,10pt]{amsart}

\usepackage{amsmath}               
\usepackage{amsthm}                
\usepackage{times}
\usepackage[all]{xy}

\numberwithin{equation}{section}

\title[Index Theory for a Class of Strongly Indefinite Functionals]%
{An Index Theory for Paths that are Solutions of a Class of Strongly Indefinite
Variational Problems.}
\author[P.\ Piccione]{Paolo Piccione}
\author[D.\ Tausk]{Daniel V.\ Tausk}
\address{Departamento de Matem\'atica,\hfill\break\indent  Universidade de S\~ao Paulo, Brazil}
\email{piccione@ime.usp.br, tausk@ime.usp.br}
\urladdr{http://www.ime.usp.br/\~{}piccione, http://www.ime.usp.br/\~{}tausk}
\thanks{The first author is partially sponsored by CNPq
(Processo n.\ 301410/95). Parts of this work were done during the visit of the two authors to the IMPA, Instituto de Matem\'atica Pura e Aplicada,
Rio de Janeiro, Brazil, in January and February 2001. The authors wish to express their
gratitude to all Faculty and Staff of the IMPA for their kind hospitality.}

\date{April 2001}

\begin{document}


\theoremstyle{plain}\newtheorem{teo}{Theorem}[section]
\theoremstyle{plain}\newtheorem{prop}[teo]{Proposition}
\theoremstyle{plain}\newtheorem{lem}[teo]{Lemma}
\theoremstyle{plain}\newtheorem{cor}[teo]{Corollary}
\theoremstyle{definition}\newtheorem{defin}[teo]{Definition}
\theoremstyle{remark}\newtheorem{rem}[teo]{Remark}
\theoremstyle{definition}\newtheorem{example}[teo]{Example}


\begin{abstract}
We generalize the  Morse index theorem of \cite{london, CRAS} and we apply the new result to obtain
lower estimates on the number of geodesics joining two fixed non conjugate points in certain classes of semi-Riemannian
manifolds. More specifically, we consider semi-Riemannian manifolds $(M,\gfrak)$ admitting a smooth distribution spanned by
commuting Killing vector fields and containing a maximal negative distribution for $\gfrak$. In particular we obtain Morse
relations for {\em stationary\/} semi-Riemannian manifolds (see \cite{GPS-JMAA}) and for the {\em G\"odel-type\/} manifolds
(see
\cite{CandelaSanchez}).
\end{abstract}

\maketitle

\begin{section}{Introduction}\label{sec:intro}
The standard Morse theory for variational problems on Hilbert manifolds
is applicable under the assumptions that the functional involved satisfies
good compactness properties (like the Palais--Smale condition), that it
has only nondegenerate critical points, and that each critical point has finite
Morse index. However, many  variational problems that arise naturally in several
contexts lead to the study of functionals that do not satisfy such assumptions.
The main example that we will keep in mind as the prototype of the theory developed
in this paper,  is the geodesic action functional on a semi-Riemannian manifold $(M,\gfrak)$,
i.e., a manifold endowed with a non positive definite nondegenerate metric tensor $\gfrak$.
It is well known that the semi-Riemannian geodesic action functional, defined
in the Hilbert manifold $\Omega_{pq}(M)$ of all curves of class $H^1$ joining
the points $p,q\in M$, does not  satisfy the Palais--Smale condition, it is unbounded
both from above and from below, and it is strongly indefinite, i.e., all its critical
points have infinite Morse index. Recall that the second variation of the
geodesic action functional at a given geodesic $\gamma$ is given by the so-called
{\em index form\/} $I_\gamma$,  which is a symmetric bounded bilinear form in the
space of all variational vector fields along $\gamma$ vanishing at the
endpoints and whose kernel consists of Jacobi fields
along $\gamma$ vanishing at the endpoints.

In order to establish the existence of at least one solution
to this kind of variational problems by means of variational methods,
a certain strategy has been suggested in some recent
works (see for instance
\cite{CandelaSanchez, GP-CAG, GPS-JMAA, Ma} and the references therein).
The basic idea in these theories is that, when the metric tensor $\gfrak$ admits
a suitable number of symmetries (i.e., Killing vector fields), then the variational
problem can be {\em reduced\/} to the study of a functional
which is bounded from below and it satisfies the Palais--Smale condition
by ``factoring out'' the negative contribution in the directions of the Killing fields.
For instance (see \cite{GP-CAG}), given a Lorentzian manifold $(M,\gfrak)$ (i.e., the index of $\gfrak$ is $1$)
admitting a timelike Killing vector field $\Ycal$, then the mentioned reduction of the variational
problem  consists in restricting the geodesic action functional to the
set of curves $\gamma\in\Omega_{pq}(M)$ that are {\em geodesics in the direction
of $\Ycal$}, i.e., that satisfy the conservation law $\gfrak(\gamma',\Ycal)\equiv C_\gamma \text{\ (constant)}$.
Similarly, for semi-Riemannian metrics $\gfrak$ of arbitrary index, the reduction is done
under the assumption that $\gfrak$ admits a family $(\Ycal_i)_{i=1}^r$ of commuting Killing vector fields
that generate a distribution which contains a maximal negative distribution
for $\gfrak$. In Ref.~\cite{GPS-JMAA} it is studied the case of stationary semi-Riemannian
manifolds, i.e., those manifolds that admit a family  of commuting Killing vector fields $\Ycal_i$
that span a maximal negative distribution for $\gfrak$, while in \cite{CandelaSanchez} it
is studied the case of Lorentzian space-times of  G\"odel type, which are Lorentzian
manifolds admitting a pair of Killing fields, whose causal character is not necessarily
constant, and that span a two-dimensional distribution where the restriction
of $\gfrak$ has index equal to $1$.

It is not hard to prove that the above mentioned {\em restricted\/} variational
problems lead to functionals whose second variation at each critical point
$\gamma$ (typically, a restriction of the index form $I_\gamma$)
is represented by a self-adjoint operator which is a compact perturbation
of a positive isomorphism, hence it has finite index. However, determining the
value of such index, as well as the problem of relating this number to the
geometrical properties of the corresponding geodesic, has turned out to be
quite a complicated task. Observe indeed that the classical Morse index theorem
fails to make sense in this situation\footnote{%
It is well known that conjugate points along a semi-Riemannian geodesic
do not in general form a finite set (see \cite{Hel1, CAG}), and the differential operator
corresponding to the Jacobi equation along the geodesic is in general
not self-adjoint.},
nor its classical proof can be adapted to obtain an alternative  statement.
In Reference~\cite{asian} the authors have proven that, in the case of
a stationary Lorentzian manifold, the index of the restricted functional
at each geodesic equals the {\em Maslov index\/} of the geodesic, which
is a homological invariant defined as an intersection number in the Grassmannian
of all Lagrangian subspaces of a symplectic space. The notion of Maslov index
was first introduced by the Russian school (see \cite{Ar} and the references
therein) and successively extended by a large number of authors in several
contexts, mainly in connection with solutions (periodic or not) of Hamiltonian
systems. In the context of semi-Riemannian geodesics, the notion of Maslov index
was introduced by Helfer in \cite{Hel1}; under generic circumstances,
the Maslov index of a semi-Riemannian geodesic is given by a sort of algebraic count
of the conjugate points (see \cite{Hel1, MPT}). Using this result,  in \cite{asian}
the authors proved the
Morse relations for geodesics of any causal character in a stationary Lorentzian
manifold, that give a lower estimate of the number of geodesics joining two
non conjugate points with a given Maslov index in terms of a Betti number
of the loop space of the underlying manifold.
The index theorem of \cite{asian} was generalized in references \cite{CRAS, semiRiem};
given a geodesic $\gamma$ in a semi-Riemannian
manifold $(M,\gfrak)$ with metric $\gfrak$ of arbitrary index, and given a maximal
negative distribution $\mathcal D_t\subset T_{\gamma(t)}M$, one defines
spaces $\mathcal K_{\mathcal D} $ and $\mathcal S_{\mathcal D}$ of variational vector fields
along $\gamma$, where $\mathcal S_{\mathcal D}$ consists of vector fields along
$\gamma$ taking values in ${\mathcal D}$ and $\mathcal K_{\mathcal D}$ consists of variational
vector fields corresponding to variations of $\gamma$ by geodesics in the directions
of ${\mathcal D}$. Then, the index $n_-\big(I_\gamma\vert_{\mathcal K_{\mathcal D}}\big)$ of the restriction of
$I_\gamma$ to $\mathcal K_{\mathcal D}$ and the
the co-index $n_+\big(I_\gamma\vert_{\mathcal S_{\mathcal D}}\big)=n_-\big({-I_\gamma\vert_{\mathcal S_{\mathcal D}}}\big)$
are finite, and their difference equals the Maslov index of $\gamma$.
In order to use this result to develop an infinite dimensional Morse theory
for semi-Riemannian geodesics using the reduction argument mentioned above, one needs to restrict to the case that
the  distribution is spanned by commuting Killing vector fields (see Subsection~\ref{sub:stationary}). Namely, in this
case the space $\mathcal K_{\mathcal D}$ is the tangent space of a Hilbert submanifold
of $\Omega_{pq}(M)$; moreover,   one has
$n_+\big(I_\gamma\vert_{\mathcal S_{\mathcal D}}\big)=0$.

The main goal of this paper is to push beyond the stationarity assumption
the results of Morse theory for semi-Riemannian geodesics, or, more
generally, for solutions of unidimensional strongly indefinite variational
problems that admit a sufficiently large number of symmetries. More precisely,
we will consider the geodesic variational problem in a semi-Riemannian
manifold $(M,\gfrak)$, with $n_-(\gfrak)=k$, that admits  linearly
independent commuting Killing vector fields $(\Ycal_i)_{i=1}^r$, with $r\ge k$,
that span a distribution ${\mathcal D}$ that {\em contains\/} a maximal
negative distribution $\Delta$. Observe that the distribution $\Delta$ itself
need not be spanned by Killing vector fields, i.e., we consider the case
that the causal character of each one of the $\Ycal_i$'s need not
be constant on $M$. Such generalization enlarges significantly the class
of variational problems to which the theory can be applied; as an example,
in this paper we will consider the case of geodesics in semi-Riemannian
manifolds of G\"odel type (Subsection~\ref{sub:Godel}).

Let us give a short description of the results of the paper, as follows.
Our main result concerning the index theory is stated and proved in the
context of {\em symplectic differential systems}. These are linear homogeneous
systems of ODE's in $\R^n\oplus{\R^n}^*$ whose coefficient matrix is
a smooth curve in the symplectic Lie algebra $\spl(2n,\R)$. Symplectic
differential systems arise naturally in connection with solutions of
Hamiltonian problems in symplectic manifolds as linearizations of
the Hamilton equations, using a symplectic referential along
the solution (see \cite[Section 3]{london}). For instance, in the geodesic case one obtains
a symplectic system (more specifically, a {\em Morse--Sturm\/} system)
by considering the Jacobi equation along the given geodesic
and using a parallel trivialization of the tangent bundle of the
semi-Riemannian manifold along the geodesic.
Associated to a symplectic system with Lagrangian initial data there
is a notion of Maslov index and of index form; in the case of
symplectic systems arising from solutions of hyper-regular Hamiltonians
(for instance, geodesics), the index form coincides with the
second variation of the Lagrangian action functional.
The choice of a distribution $\Dcal$ determines a {\em reduced\/}
symplectic system; we prove that when $\Dcal$ contains a maximal
negative distribution then  the Morse index of the restriction
of the index form to the space $\Kcal_{\Dcal}$ equals the difference
between the Maslov index of the original symplectic system and
the Maslov index of the reduced symplectic system, plus a suitable
correction term determined by the initial condition (Theorem~\ref{thm:INDEXTHMNOVO}).

The paper is organized as follows. In Section~\ref{sec:sds} we recall from
\cite{london} the general formalism of symplectic differential systems and
the index theorem of \cite{CRAS}. In Section~\ref{sec:IndexTheorem} we state
and prove our generalized index theorem; finally, in Section~\ref{sec:semiRiem} we
discuss the applications to semi-Riemannian geodesics.
\end{section}

\begin{section}{Symplectic Differential Systems}\label{sec:sds}

We start by introducing the basic notation and terminology that
will be used throughout the paper. Given real vector spaces $V$,
$W$ we denote by $\Lin(V,W)$ the space of linear maps from $V$ to
$W$ and by $\Bil(V,W)$ the space of bilinear forms $B:V\times
W\to\R$; by $\Bils(V)$ we denote the subspace of $\Bil(V,V)$
consisting of {\em symmetric\/} bilinear forms. For
$T\in\Lin(V,W)$ we denote by $T^*\in\Lin(W^*,V^*)$ the {\em
transpose\/} of $T$, where $V^*$, $W^*$ denote respectively the
dual spaces of $V$ and $W$. The {\em index\/} of a symmetric
bilinear form $B\in\Bils(V)$, denoted by $n_-(B)$, is defined as
the supremum of the dimensions of the subspaces of $V$ on which
$B$ is negative definite:
\[n_-(B)=\sup\big\{\Dim(W):B\vert_W\ \text{is negative definite}\big\}\in\N\cup\{+\infty\};\]
the {\em coindex\/} of $B$ is defined by $n_+(B)=n_-(-B)$ and the
{\em signature\/} of $B$ is defined as the difference
$\sgn(B)=n_+(B)-n_-(B)$, provided that one of the numbers
$n_+(B)$, $n_-(B)$ is finite. We also define the {\em
degeneracy\/} of a symmetric bilinear form $B$ as the dimension of
its kernel, i.e., $\dgn(B)=\Dim\big(\Ker(B)\big)$.

We always implicitly identify the spaces $\Bil(V,W)$ and
$\Lin(V,W^*)$ by the natural isomorphism $B(v,w)=B(v)(w)$. With
such identification, if $V$ is finite dimensional, then a bilinear
form $B\in\Bil(V,V)\cong\Lin(V,V^*)$ is symmetric iff $B$ equals
its own transpose
$B^*\in\Lin(V^{**},V^*)\cong\Lin(V,V^*)\cong\Bil(V,V)$; moreover,
if $B\in\Bils(V)$ is {\em nondegenerate\/} (i.e., $\Ker(B)=\{0\}$)
then $B^{-1}\in\Lin(V^*,V)\cong\Lin(V^*,V^{**})\cong\Bil(V^*,V^*)$
is the nondegenerate symmetric bilinear form on $V^*$ which equals
the push-forward of $B\in\Bils(V)$ by the isomorphism $B:V\to
V^*$.

A central object for our theory is the symplectic space $\R^n\oplus{\R^n}^*$ endowed with the canonical symplectic form
\[\omega\big((v_1,\alpha_1),(v_2,\alpha_2)\big)=\alpha_2(v_1)-\alpha_1(v_2);\]
we denote by $\Spl(2n,\R)$ the symplectic group of
$(\R^n\oplus{\R^n}^*,\omega)$, which is the closed subgroup of the
general linear group $\GL(2n,\R)$ consisting of those isomorphisms
of $\R^n\oplus{\R^n}^*$ that preserve $\omega$. The Lie algebra
$\spl(2n,\R)$ of $\Spl(2n,\R)$ consists of the linear
endomorphisms $X$ of $\R^n\oplus{\R^n}^*$ that are represented by
matrices of the form:
\begin{equation}\label{eq:XABC}
X=\begin{pmatrix}A&B\\C&-A^*\end{pmatrix},
\end{equation}
with $B$, $C$ symmetric. We think of $A$ as a linear endomorphism
of $\R^n$ and of $B$, $C$ as linear maps $B:{\R^n}^*\to\R^n$,
$C:\R^n\to{\R^n}^*$; we also think of $B$ as a symmetric bilinear
form on ${\R^n}^*$ and $C$ as a symmetric bilinear form on $\R^n$.
We can now give the following:
\begin{defin}\label{thm:defsds}
Let $X:[a,b]\to\spl(2n,\R)$ be a smooth curve and define $A$, $B$ and $C$ as in \eqref{eq:XABC}. Assume that $B(t)$ is
nondegenerate for all $t\in[a,b]$. The {\em symplectic differential system\/} in $\R^n$ with {\em coefficient matrix\/}
$X$ is the following linear homogeneous first order system of ODE's:
\begin{equation}\label{eq:sds}
\frac{\dd}{\dd
t}\begin{pmatrix}v(t)\\\alpha(t)\end{pmatrix}=X(t)\begin{pmatrix}v(t)\\\alpha(t)\end{pmatrix}.
\end{equation}
\end{defin}
We identify the symplectic differential system \eqref{eq:sds} with
its coefficient matrix $X$ so that we will in general refer to
{\em the symplectic differential system $X$}; the blocks $A$, $B$,
$C$ will be called the {\em coefficients\/} of $X$. The
nondegeneracy of $B(t)$ implies that the index of $B(t)$ is
independent of $t\in[a,b]$; we call $B$ the {\em fundamental
coefficient\/} of the symplectic differential system $X$ and the
integer $n_-\big(B(t)\big)$ the {\em index\/} of $X$.

\begin{example}\label{thm:exaMS}
A special class of symplectic differential systems are the so called {\em Morse--Sturm\/} systems; in our notation,
Morse--Sturm systems are symplectic differential systems with $A\equiv0$. They can be written in the form of a second order
linear differential equation:
\begin{equation}\label{eq:MS}
g^{-1}\big(gv'\big)'=Rv,
\end{equation}
where $g:[a,b]\to\Bils(\R^n)$, $R:[a,b]\to\Lin(\R^n)$ are smooth curves, $g(t)$ is nondegenerate and $g(t)R(t)$ is symmetric
for all $t\in[a,b]$. Equation \eqref{eq:MS} is identified with the symplectic differential system with coefficients $A=0$,
$B=g^{-1}$ and $C=gR$. Observe that when $g$ is constant, \eqref{eq:MS} becomes:
\begin{equation}\label{eq:MSJacobi}
v''=Rv.
\end{equation}
\end{example}

Let $X$ be a symplectic differential system. Given a map $v:[a,b]\to\R^n$ there exists at most one map
$\alpha:[a,b]\to{\R^n}^*$ for which
$(v,\alpha)$ can be a solution of $X$; such map will be denoted by $\alpha_v$ and it is given by:
\[\alpha_v=B^{-1}(v'-Av).\]
We will usually say that a map $v:[a,b]\to\R^n$ is a solution of $X$ meaning that $(v,\alpha_v)$ is a solution of $X$. For
later use, we give here the following formula:
\begin{equation}\label{eq:alphafv}
\alpha_{fv}=f\alpha_v+f'B^{-1}v,
\end{equation}
valid for all absolutely continuous maps $v:[a,b]\to\R^n$, $f:[a,b]\to\R$.

Recall that a {\em Lagrangian subspace\/} of $\R^n\oplus{\R^n}^*$
is an $n$-dimensional subspace on which the symplectic form
$\omega$ vanishes. We will {\em henceforth\/} denote by $L_0$ the
Lagrangian subspace:
\begin{equation}\label{eq:defL0}
L_0=\{0\}\oplus{\R^n}^*\subset\R^n\oplus{\R^n}^*.
\end{equation}
\begin{defin}\label{thm:definitial}
Given a Lagrangian
subspace
$\ell_0\subset\R^n\oplus{\R^n}^*$ we consider the following initial condition for the system \eqref{eq:sds}:
\begin{equation}\label{eq:ell0IC}
\big(v(a),\alpha(a)\big)\in\ell_0;
\end{equation}
the pair $(X,\ell_0)$ will be called a {\em symplectic
differential system with initial data}. We call a pair
$(v,\alpha)$ an {\em $(X,\ell_0)$-solution\/} (or a {\em solution
of $(X,\ell_0)$}) if $(v,\alpha)$ satisfies \eqref{eq:sds} and
\eqref{eq:ell0IC}.
\end{defin}
It is easy to see that a Lagrangian
$\ell_0\subset\R^n\oplus{\R^n}^*$ determines uniquely a pair
$(P,S)$ where $P\subset\R^n$ is a subspace, $S:P\to P^*$ is a
symmetric bilinear form on $P$ and
\begin{equation}\label{eq:ell0PS}
\ell_0=\big\{(v,\alpha):v\in P,\ \alpha\vert_P+S(v)=0\in P^*\big\};
\end{equation}
in terms of the pair $(P,S)$ the initial condition
\eqref{eq:ell0IC} can be rewritten as:
\begin{equation}\label{eq:ICPS}
v(a)\in P,\quad\alpha(a)\vert_P+S\big(v(a)\big)=0\in P^*.
\end{equation}
The initial condition \eqref{eq:ell0IC} (or \eqref{eq:ICPS}) is
called {\em nondegenerate\/} if $B(a)^{-1}\in\Bils(\R^n)$ is
nondegenerate on $P$. This is always the case, for instance, if
$\ell_0=L_0$; in this case $P=\{0\}$ and \eqref{eq:ICPS} reduces
to $v(a)=0$.

We denote by $\mathbb V$ the $n$-dimensional space of all solutions $v$ of $(X,\ell_0)$; for $t\in[a,b]$ we set:
\[\mathbb V[t]=\big\{v(t):v\in\mathbb V\big\}\subset\R^n.\]

For each $t\in[a,b]$, we have an isomorphism
$\Phi(t)$ of $\R^n\oplus{\R^n}^*$ defined by the relation $\Phi(t)\big(v(a),\alpha(a)\big)=\big(v(t),\alpha(t)\big)$ for
every solution $(v,\alpha)$ of $X$;   $\Phi$ satisfies the matrix differential equation
\begin{equation}\label{eq:defPhit}
\Phi'(t)=X(t)\Phi(t),
\end{equation}
with initial condition $\Phi(a)=\Id$. It follows that $\Phi$ is a smooth curve in the symplectic group
$\Spl(2n,\R)$. We call $\Phi$ the {\em fundamental matrix\/} of the system $X$. For $t\in[a,b]$, we set:
\begin{equation}\label{eq:defellt}
\ell(t)=\Phi(t)(\ell_0)=\big\{\big(v(t),\alpha_v(t)\big):v\in\mathbb V\big\};
\end{equation}
obviously $\ell(t)$ is a Lagrangian subspace of $\R^n\oplus{\R^n}^*$.

\begin{defin}\label{thm:defconjs}
A {\em focal instant\/} for the symplectic differential system
with initial data $(X,\ell_0)$ is an instant
$t\in\left]a,b\right]$ such that there exists a non zero
$(X,\ell_0)$-solution $v$ with $v(t)=0$; the dimension of the
space of all $(X,\ell_0)$-solutions vanishing at $t$ is called the
{\em multiplicity\/} of the focal instant $t$ and is denoted by
$\mul(t)$. The {\em signature\/} of a focal instant
$t\in\left]a,b\right]$, denoted by $\sgn(t)$, is defined as the
signature of the restriction of $B(t)^{-1}$ to the
$B(t)^{-1}$-orthogonal complement of $\mathbb V[t]$; the focal
instant is called nondegenerate when such restriction is
nondegenerate. In the special case where $\ell_0=L_0$ then the
focal instants of $(X,\ell_0)$ are also called {\em conjugate
instants\/} of $X$.
\end{defin}
An instant $t\in\left]a,b\right]$ is focal iff $\mathbb V[t]$ is a
proper subspace of $\R^n$, in which case $\mul(t)$ is the
codimension of $\mathbb V[t]$. It is well known (see for instance
\cite[Theorem 2.3.3]{london}) that nondegenerate focal instants
are isolated.

Formula \eqref{eq:defellt} defines a smooth curve in the {\em
Lagrangian Grassmannian\/} $\Lambda$ of the symplectic space
$\R^n\oplus{\R^n}^*$ which is the embedded submanifold of the
Grassmannian of all $n$-dimensional subspaces of
$\R^n\oplus{\R^n}^*$ consisting of all Lagrangian subspaces. We
denote by $\Lambda^{\ge1}(L_0)\subset\Lambda$ the {\em Maslov
cycle\/} of $L_0$ which is the set of all Lagrangians
$L\in\Lambda$ with $L\cap L_0\ne\{0\}$. Each continuous curve $l$
in $\Lambda$ with endpoints outside $\Lambda^{\ge1}(L_0)$ defines
a relative singular homology class with integer coefficients in
$H_1\big(\Lambda,\Lambda\setminus\Lambda^{\ge1}(L_0)\big)\cong\Z$
and the corresponding integer number $\mu_{L_0}(l)$ is called the
{\em Maslov index\/} of the curve $l$. If $\ell_0$ defines a
nondegenerate initial condition for $X$ then there are no
$(X,\ell_0)$-focal instants near $t=a$ and we can give the
following:
\begin{defin}\label{thm:defMaslov}
Let $(X,\ell_0)$ be a symplectic differential system with initial data; assume that $t=b$ is not focal and that $\ell_0$
defines a nondegenerate initial condition for $X$. The {\em Maslov index\/} of $(X,\ell_0)$ is defined by:
\[\iMaslov(X,\ell_0)=\mu_{L_0}\big(\ell\vert_{[a+\varepsilon,b]}\big),\]
where $\varepsilon>0$ is chosen in such a way that there are no
focal instants in $\left]a,a+\varepsilon\right]$. When
$\ell_0=L_0$ we call $\iMaslov(X,\ell_0)$ the {\em Maslov index of
$X$\/} and we write just $\iMaslov(X)$ .
\end{defin}
For more details concerning the geometry of the Lagrangian
Grassmannian $\Lambda$ of a symplectic space and the Maslov index
of curves in $\Lambda$ we refer to \cite{Duis, MPT, RobSal}; for more
details concerning the notion of Maslov index for symplectic
differential systems with initial data we refer to
\cite[Subsection 2.3]{london}. The Maslov index of a pair
$(X,\ell_0)$ is stable by uniformly small perturbations of $X$.
Generically, it gives an algebraic count of the focal instants of
$(X,\ell_0)$:
\begin{prop}\label{thm:calcMaslov}
If all the focal instants of $(X,\ell_0)$ are nondegenerate, $t=b$ is not focal and $\ell_0$ defines a nondegenerate initial
condition for $X$ then:
\[\iMaslov(X,\ell_0)=\sum_{t\in\left]a,b\right[}\sgn(t).\]
\end{prop}
\begin{proof}
See \cite[Theorem 2.3.3]{london}.
\end{proof}
The nondegeneracy assumption on the focal instants of $(X,\ell_0)$
is indeed generic (see \cite[Proposition 2.4.1]{london}); in the
case of degenerate focal instants, the Maslov index can also be
explicitly computed in terms of suitable coordinate charts on the
Lagrangian Grassmannian (see \cite[Proposition 4.3.1]{MPT}).

Given a symplectic differential system with initial data $(X,\ell_0)$ we consider the Hilbert space:
\[\Hcal=\big\{v\in H^1([a,b],\R^n):v(a)\in P,\ v(b)=0\big\},\]
where $H^1([a,b],\R^n)$ denotes the Sobolev space of absolutely
continuous $\R^n$-valued maps on $[a,b]$ with square integrable
derivatives. On $\Hcal$ we define a bounded bilinear form $I$
given by:
\begin{equation}\label{eq:indexform}
I(v,w)=\int_a^bB(\alpha_v,\alpha_w)+C(v,w)\,\dd
t-S\big(v(a),w(a)\big).
\end{equation}
We call $I$ the {\em index form\/} associated to the pair
$(X,\ell_0)$. The index form of $(X,L_0)$ will be simply referred
to as the {\em index form of $X$}; the term $S\big(v(a),w(a)\big)$
in \eqref{eq:indexform} does not appear in the index form of $X$
and the boundary conditions defining $\Hcal$ become $v(a)=v(b)=0$.

We recall from \cite[Subsection 2.10]{london} the notion of
isomorphism in the class of symplectic differential systems:
\begin{defin}\label{thm:defiso}
Let $X$, $\tilde X$ be symplectic differential systems. An {\em
isomorphism\/} $\phi:X\to\tilde X$ is a smooth curve
$\phi:[a,b]\to\Spl(2n,\R)$ with $\phi(t)(L_0)=L_0$ for all $t$ and
such that one of the following equivalent conditions are
satisfied:
\begin{itemize}
\item[(a)] $\tilde X(t)=\phi'(t)\phi(t)^{-1}+\phi(t)X(t)\phi(t)^{-1}$ for all $t\in[a,b]$;

\item[(b)] $\tilde\Phi(t)=\phi(t)\Phi(t)\phi(a)^{-1}$ for all $t$, where $\Phi$, $\tilde\Phi$ denote
respectively the fundamental matrices of $X$ and $\tilde X$;

\item[(c)] $(v,\alpha)$ is a solution of $X$ iff $(\tilde v,\tilde\alpha)=\phi(v,\alpha)$ is a solution of
$\tilde X$.

\end{itemize}
Given Lagrangians $\ell_0,\tilde\ell_0\subset\R^n\oplus{\R^n}^*$,
we say that $\phi:(X,\ell_0)\to(\tilde X,\tilde\ell_0)$ is an {\em
isomorphism\/} if $\phi:X\to\tilde X$ is an isomorphism and
$\phi(a)(\ell_0)=\tilde\ell_0$.
\end{defin}
Elements $\phi\in\Spl(2n,\R)$ that preserve $L_0$ are written in block matrix form as:
\begin{equation}\label{eq:phiZW}
\phi=\begin{pmatrix}Z&0\\{Z^*}^{-1}W&{Z^*}^{-1}\end{pmatrix},
\end{equation}
with $Z\in\Lin(\R^n,\R^n)$ invertible and $W\in\Bils(\R^n)$.

\begin{prop}\label{thm:tudoinvariante}
Let $\phi:(X,\ell_0)\to(\tilde X,\tilde\ell_0)$ be an isomorphism and define $Z$, $W$ as in \eqref{eq:phiZW}. Then:
\begin{itemize}
\item $(X,\ell_0)$ and $(\tilde X,\tilde\ell_0)$ have the same focal instants with same multiplicity and signature;

\item a focal instant is nondegenerate for $(X,\ell_0)$ iff it is nondegenerate for $(\tilde X,\tilde\ell_0)$;

\item $\ell_0$ defines a nondegenerate initial condition for $X$ iff $\tilde\ell_0$ defines a nondegenerate initial condition
for $\tilde X$;

\item if $t=b$ is not focal and $\ell_0$ defines a nondegenerate initial condition for $X$ then
$\iMaslov(X,\ell_0)=\iMaslov(\tilde X,\tilde\ell_0)$;

\item the map $\Hcal\ni v\mapsto Zv\in\tilde{\Hcal}$ is a continuous isomorphism and it carries the index form $I:\Hcal\times\Hcal\to\R$ of $(X,\ell_0)$ to the
index form $\tilde I:\tilde\Hcal\times\tilde\Hcal\to\R$ of
$(\tilde X,\tilde\ell_0)$.
\end{itemize}
\end{prop}
\begin{proof}
See \cite[Subsection 2.10]{london}.
\end{proof}

\begin{defin}\label{thm:distrib}
For each $t\in[a,b]$, let $\Dcal_t$ be an $r$-dimensional subspace
of $\R^n$, where $r=0,\ldots,n$ is fixed. We call $\Dcal$ a {\em
smooth family of subspaces\/} in $\R^n$ over the interval $[a,b]$
if there exist smooth maps $Y_i:[a,b]\to\R^n$, $i=1,\ldots,r$,
such that $\big(Y_i(t)\big)_{i=1}^r$ is a basis for $\Dcal_t$ for
every $t\in[a,b]$; the integer $r$ is called the {\em rank\/} of
the family $\Dcal$ and the family of maps $(Y_i)_{i=1}^r$ is
called a {\em frame\/} for $\Dcal$. A smooth family of subspaces
$\Dcal$ is called {\em nondegenerate\/} (resp., {\em negative})
for a symplectic differential system $X$ if the symmetric bilinear
form $B(t)^{-1}\in\Bils(\R^n)$ is nondegenerate (resp., negative
definite) on $\Dcal_t$ for all $t\in[a,b]$. If $\Dcal$ is negative
for $X$ and the rank of $\Dcal$ equals the index of $X$ we say
that $\Dcal$ is {\em maximal negative\/} for $X$. If $\Dcal$ is
nondegenerate for $X$ then the {\em index\/} of $\Dcal$ with
respect to $X$ is defined to be the index of the restriction of
$B(t)^{-1}$ to $\Dcal_t$ (which does not depend on $t\in[a,b]$).
\end{defin}

Let $\Dcal$ be a fixed smooth family of subspaces of rank $r$ in
$\R^n$ over the interval $[a,b]$. A map $Y:[a,b]\to\R^n$ is called
a {\em section\/} of $\Dcal$ if $Y(t)\in\Dcal_t$ for all
$t\in[a,b]$; we set:
\begin{equation}\label{eq:SD}
\Scal_\Dcal=\big\{v\in\Hcal:\text{$v$ is a section of $\Dcal$ and $v(a)=0$}\big\}.
\end{equation}
An absolutely
continuous map
$v:[a,b]\to\R^n$ is called a {\em solution of
$X$ along
$\Dcal$\/} if for every absolutely continuous section $Y:[a,b]\to\R^n$ of $\Dcal$ we have that $t\mapsto\alpha_v(t)Y(t)$ is
absolutely continuous and that:
\begin{equation}\label{eq:alongD}
\big(\alpha_v(Y)\big)'=B(\alpha_v,\alpha_Y)+C(v,Y);
\end{equation}
we set:
\begin{equation}\label{eq:KD}
\Kcal_\Dcal=\big\{v\in\Hcal:\text{$v$ is a solution of $X$ along $\Dcal$}\big\}.
\end{equation}
When $v'$ is also absolutely continuous then the left hand side of
\eqref{eq:alongD} can be expanded and one obtains that $v$ is a
solution of $X$ along $\Dcal$ iff $\alpha_v'-Cv+A^*\alpha_v$
vanishes on $\Dcal$; in particular, every solution $v$ of $X$ is a
solution along $\Dcal$. If $(Y_i)_{i=1}^r$ is a frame for $\Dcal$
then $v$ is a solution of $X$ along $\Dcal$ iff $\alpha_v(Y_i)$ is
absolutely continuous and \eqref{eq:alongD} is satisfied with
$Y=Y_i$ for all $i=1,\ldots,r$.

The spaces $\Kcal_\Dcal$ and $\Scal_\Dcal$ are orthogonal with
respect to the index form:
\begin{lem}\label{thm:KdSdorthog}
Consider a triple $(X,\ell_0,\Dcal)$ where $(X,\ell_0)$ is a
symplectic differential system with initial data and $\Dcal$ is a
smooth family of subspaces in $\R^n$ over the interval $[a,b]$.
Then $I(v,w)=0$ for all $v\in\Kcal_\Dcal$ and all
$w\in\Scal_\Dcal$.
\end{lem}
\begin{proof}
Choose a frame $(Y_i)_{i=1}^r$ for $\Dcal$ and write
$w=\sum_{i=1}^rf_iY_i$; the conclusion follows from an easy
computation using \eqref{eq:alphafv}, \eqref{eq:indexform} and
\eqref{eq:alongD}.
\end{proof}

\begin{defin}\label{thm:reducedsystem}
Given a nondegenerate smooth family of subspaces $\Dcal$ for a
symplectic differential system $X$ and given a frame
$(Y_i)_{i=1}^r$ for $\Dcal$ we define the {\em reduced symplectic
differential system\/} corresponding to $X$, $\Dcal$ and
$(Y_i)_{i=1}^r$ to be the following symplectic differential system
in $\R^r$:
\begin{equation}\label{eq:sistred}
\left\{
\begin{aligned}
f'&=-(\mathfrak B^{-1}\circ\mathcal A)f+\mathfrak B^{-1}\varphi,\\
\varphi'&=(\mathcal C-\mathcal A^*\circ\mathfrak B^{-1}\circ\mathcal A)f+(\mathcal A^*\circ\mathfrak B^{-1})\varphi,
\end{aligned}
\right.
\end{equation}
where $\mathcal A(t),\mathfrak B(t),\mathcal C(t)\in\Lin(\R^r,{\R^r}^*)$ are the linear operators represented by the following
matrices:
\begin{equation}\label{eq:calABC}
\begin{split}
\mathfrak B_{ij}&(t)=B(t)^{-1}\big(Y_i(t),Y_j(t)\big),\qquad\mathcal A_{ij}(t)=\alpha_{Y_j}(t)Y_i(t)\\
&\!\!\mathcal C_{ij}(t)=B(t)\big(\alpha_{Y_i}(t),\alpha_{Y_j}(t)\big)+C(t)\big(Y_i(t),Y_j(t)\big).
\end{split}
\end{equation}
\end{defin}
The fact that $\Dcal$ is nondegenerate for $X$ implies that $\mathfrak B(t)$ is indeed nondegenerate, so that $\mathfrak
B(t)^{-1}$ in \eqref{eq:sistred} makes sense. Moreover, the index of the reduced system \eqref{eq:sistred} equals the index of
$\Dcal$ with respect to $X$. We will denote by $X_\red$ the coefficient matrix of \eqref{eq:sistred} and by $A_\red$,
$B_\red$, $C_\red$ the coefficients of $X_\red$.

The introduction of the reduced symplectic system is motivated by the following:
\begin{lem}\label{thm:whyreduced}
If $\Dcal$ is a nondegenerate smooth family of subspaces for $X$
and $(Y_i)_{i=1}^r$ is a frame for $\Dcal$ then
$v=\sum_{i=1}^rf_iY_i$ is a solution of $X$ along $\Dcal$ iff
$f=(f_i)_{i=1}^r$ is a solution of the reduced symplectic system
$X_\red$ corresponding to $\Dcal$ and $(Y_i)_{i=1}^r$. In
particular, $\Kcal_\Dcal\cap\Scal_\Dcal=\{0\}$ iff $t=b$ is not
conjugate for $X_\red$.
\end{lem}
\begin{proof}
It is a straightforward computation using \eqref{eq:alphafv}.
\end{proof}

The index form of the reduced symplectic system can be identified with a restriction of the index form of the original system:
\begin{lem}\label{thm:identifyindexformred}
If $I_\red:\Hcal_\red\times\Hcal_\red\to\R$ denotes the index form
of $X_\red$ then the continuous isomorphism
\begin{equation}\label{eq:isolambda}
\lambda:\Hcal_\red\ni
f=(f_i)_{i=1}^r\longmapsto\sum_{i=1}^rf_iY_i\in\Scal_\Dcal
\end{equation}
carries $I_\red$ to $I\vert_{\Scal_\Dcal}$, i.e.,
$I\big(\lambda(f),\lambda(g)\big)=I_\red(f,g)$ for all
$f,g\in\Hcal_\red$.
\end{lem}
\begin{proof}
Recall that the domain $\Hcal_\red$ of the index form of $X_\red$
is given by:
\[\Hcal_\red=\big\{f\in H^1([a,b],\R^r):f(a)=f(b)=0\big\},\]
so that $\lambda$ is indeed a continuous isomorphism. The
conclusion follows by a straightforward computation using
\eqref{eq:alphafv}.
\end{proof}

\begin{rem}\label{thm:mudaframe}
The reduced symplectic system $X_\red$ depends on the choice of the frame for $\Dcal$ and not only on $\Dcal$; however,
different choices of a frame for $\Dcal$ produce isomorphic reduced systems. In particular, when discussing notions that are
invariant by isomorphisms (like Maslov index, focal instants), we do not need to specify a frame for $\Dcal$.
\end{rem}


We recall the following result announced in \cite{CRAS} and proven
in \cite{semiRiem}:
\begin{teo}\label{thm:INDEXvelho}
Consider a triple $(X,\ell_0,\Dcal)$, where $(X,\ell_0)$ is a
symplectic differential system with initial data such that
$\ell_0$ defines a nondegenerate initial condition for $X$ and
$\Dcal$ is a maximal negative smooth family of subspaces for $X$.
If $t=b$ is not $(X,\ell_0)$-focal then:
\[\iMaslov(X,\ell_0)=n_-\big(I\vert_{\Kcal_\Dcal}\big)-n_+\big(I\vert_{\Scal_\Dcal}\big)-n_-\big(B(a)^{-1}\vert_P\big).\]
\qed
\end{teo}

The following lemma\footnote{%
In \cite[Corollary 2.6.10]{london} the result of
Lemma~\ref{thm:HKmaisS} is proven for a maximal negative family
$\Dcal$; here we adapt the proof to the case where $\Dcal$ is only
nondegenerate.} is an addendum to the result of
Theorem~\ref{thm:INDEXvelho}:
\begin{lem}\label{thm:HKmaisS}
Let $(X,\ell_0,\Dcal)$ be a triple where $(X,\ell_0)$ is a
symplectic differential system with initial data and $\Dcal$ is a
nondegenerate smooth family of subspaces for $X$. If $t=b$ is not
a conjugate instant for the reduced symplectic system
corresponding to $(X,\ell_0)$ and $\Dcal$ then we have a direct
sum decomposition $\Hcal=\Kcal_\Dcal\oplus\Scal_\Dcal$.
\end{lem}
\begin{proof}
If $(Y_i)_{i=1}^r$ denotes a frame for $\Dcal$ then $\Kcal_\Dcal$
is the kernel of the bounded linear operator $F:\Hcal\to
L^2\big([a,b],{\R^r}^*\big)/\cte$ defined by
\[F(v)(t)_i=\alpha_v(t)Y_i(t)-\int_a^tB(\alpha_v,\alpha_{Y_i})+C(v,Y_i)\,\dd
s\mod{\cte},\] for all $v\in\Hcal$, $t\in[a,b]$, $i=1,\ldots,r$,
where $L^2\big([a,b],{\R^r}^*\big)$ denotes the Hilbert space of
${\R^r}^*$-valued square integrable maps on $[a,b]$ and $\cte$
denotes the $r$-dimensional subspace of
$L^2\big([a,b],{\R^r}^*\big)$ consisting of constant maps. If
$\lambda$ denotes the isomorphism defined in \eqref{eq:isolambda}
then a straightforward computation using \eqref{eq:alphafv} shows
that $F\circ\lambda$ is given by:
\begin{equation}\label{eq:Fbolalambda}
(F\circ\lambda)(f)(t)=\mathfrak B(t)f'(t)+\mathcal
A(t)f(t)-\int_a^t\mathcal A^*f'+\mathcal Cf\,\dd s\mod{\cte},
\end{equation}
for all $f=(f_i)_{i=1}^r\in\Hcal_\red$, $t\in[a,b]$. The kernel of
$F\circ\lambda$ are the solutions $f$ of the reduced symplectic
system with $f(a)=f(b)=0$; it follows that $F\circ\lambda$ is
injective. We will now show that $F\circ\lambda$ is a Fredholm
operator of index zero; this will imply that $F\circ\lambda$ (and
therefore $F\vert_{\Scal_\Dcal}$) is an isomorphism and that will
conclude the proof. To prove that $F\circ\lambda$ is a Fredholm
operator of index zero, observe first that by the additivity of
the Fredholm index by composition of operators and by the
invertibility of $\mathfrak B$, the map
\begin{equation}\label{eq:frakBflinha}
\Hcal_\red\ni f\longmapsto\mathfrak Bf'\mod{\cte}\in
L^2\big([a,b],{\R^r}^*\big)/\cte
\end{equation}
is a Fredholm operator of index zero. By \eqref{eq:Fbolalambda}
and the compact inclusion $H^1\hookrightarrow C^0$, the operator
$F\circ\lambda$ is a compact perturbation of
\eqref{eq:frakBflinha}. This concludes the proof.
\end{proof}

\begin{rem}\label{thm:KerI}
The kernel of the index form $I:\Hcal\times\Hcal\to\R$ is the
space of solutions $v:[a,b]\to\R^n$ of $(X,\ell_0)$ with $v(b)=0$
(see \cite[Subsection 2.5]{london}). It follows from
Lemmas~\ref{thm:KdSdorthog} and \ref{thm:HKmaisS} that if $t=b$ is
not conjugate for the reduced symplectic system then the kernel of
$I\vert_{\Kcal_\Dcal}$ coincides with the kernel of $I$ in
$\Hcal$.
\end{rem}

\begin{rem}\label{thm:outrored}
In some situations it is useful to consider the symplectic
differential system $\tilde X_\red$ which is isomorphic to
$X_\red$ and whose coefficients $\tilde A_\red$, $\tilde B_\red$,
$\tilde C_\red$ are given by:
\begin{gather}
\label{eq:tildeABCred}\tilde A_\red(t)=-\mathfrak B(t)^{-1}\circ\mathcal A_{\mathrm{ant}}(t),\quad\tilde B_\red=\mathfrak B(t)^{-1},\hphantom{\text{\eqref{eq:tildeABCred}}}\\
\notag\tilde C_\red(t)=\mathcal C(t)-\mathcal A'_{\mathrm{sym}}(t)+\mathcal A_{\mathrm{ant}}(t)\circ\mathfrak
B(t)^{-1}\circ\mathcal A_{\mathrm{ant}}(t),
\end{gather}
for all $t\in[a,b]$, where $\mathcal A_{\mathrm{sym}}$, $\mathcal
A_{\mathrm{ant}}$ denote respectively the symmetric and
anti-symmetric components of $\mathcal A$:
\begin{equation}\label{eq:Aantsim}
\mathcal A_{\mathrm{sym}}(t)=\frac{\mathcal A(t)+\mathcal
A(t)^*}2,\quad\mathcal A_{\mathrm{ant}}(t)=\frac{\mathcal A(t)-\mathcal
A(t)^*}2.
\end{equation}
An explicit isomorphism from $X_\red$  to $\tilde X_\red$  is given by \eqref{eq:phiZW} with:
\[Z(t)=\Id,\quad W(t)=-\mathcal A_{\mathrm{sym}}(t),\qquad t\in[a,b].\]
\end{rem}

\begin{rem}\label{thm:simplyred}
If a nondegenerate smooth family of subspaces $\Dcal$ for a
symplectic differential system $X$ admits a frame $(Y_i)_{i=1}^r$
consisting of {\em solutions\/} of $X$ satisfying the symmetry
condition
\[\alpha_{Y_i}(Y_j)=\alpha_{Y_j}(Y_i),\quad i,j=1,\ldots,r,\]
then the coefficients of the reduced symplectic system $\tilde
X_\red$ defined in Remark~\ref{thm:outrored} are $\tilde
A_\red=0$, $\tilde B_\red=\mathfrak B^{-1}$, $\tilde C_\red=0$.
The system $\tilde X_\red$ becomes the differential equation
\begin{equation}\label{eq:veryreduced}
\mathfrak Bf'\equiv\text{constant}.
\end{equation}
An instant $t\in\left]a,b\right]$ is conjugate for $\tilde X_\red$
iff the integral:
\begin{equation}\label{eq:Bintegral}
\mathbfB(t)=\int_a^t\mathfrak B(s)^{-1}\,\dd s
\end{equation}
is a degenerate (symmetric) bilinear form in ${\R^r}^*$, in which case the multiplicity of $t$ equals the degeneracy of
$\mathbfB(t)$. If $t=b$ is not conjugate for $\tilde X_\red$ then, by Proposition~\ref{thm:calcMaslov}, the Maslov
index of $\tilde X_\red$ is given by:
\begin{equation}\label{eq:Maslovsimplyred}
\iMaslov(\tilde
X_\red)=\sum_{t\in\left]a,b\right[}\sgn\Big(\mathfrak
B(t)\vert_{\Img(\mathbfB(t))^\perp}\Big),
\end{equation}
provided that $\mathfrak B(t)$ is nondegenerate on the image of
$\mathbfB(t)$ for those $t\in\left]a,b\right]$ such that
$\mathbfB(t)$ is degenerate. In \eqref{eq:Maslovsimplyred} we have
denoted by $\perp$ the orthogonal complement with respect to
$\mathfrak B(t)$.
\end{rem}

\end{section}

\begin{section}{The Index Theorem}\label{sec:IndexTheorem}

In this section we consider the following setup:
\begin{itemize}
\item $(X,\ell_0)$ is a symplectic differential system with initial data on $\R^n$ over the interval $[a,b]$;

\item $\Dcal$ and $\Delta$ are nondegenerate smooth families of subspaces for $X$ with $\Delta_t\subset\Dcal_t$ for all $t$;

\item $(Y_i)_{i=1}^r$ is a frame for $\Dcal$ such that $(Y_i)_{i=1}^k$ is a frame for $\Delta$;


\item $X_\red$ is the reduced symplectic system corresponding to $\Dcal$ and $(Y_i)_{i=1}^r$;

\item $I:\Hcal\times\Hcal\to\R$ is the index form of $(X,\ell_0)$ and $\Kcal_\Dcal$, $\Scal_\Dcal$ (resp., $\Kcal_\Delta$,
$\Scal_\Delta$) are the subspaces of $\Hcal$ defined in analogy with \eqref{eq:KD} and \eqref{eq:SD} for $(X,\ell_0)$ and
$\Dcal$ (resp.,
$\Delta$);


\item $I_\red:\Hcal_\red\times\Hcal_\red\to\R$ is the index form of $X_\red$;

\item $\lambda:\Hcal_\red\to\Scal_\Dcal$ is the continuous isomorphism defined in \eqref{eq:isolambda};

\item $\Delta_\red$ is the (constant) smooth family of subspaces $\Delta_\red\equiv\R^k\oplus\{0\}\subset\R^r$ in $\R^r$ over the interval $[a,b]$;

\item $\Kcal_{\Delta_\red}$ and $\Scal_{\Delta_\red}$ are the subspaces of $\Hcal_\red$ defined in analogy with \eqref{eq:KD}
and \eqref{eq:SD} for the symplectic differential system with
initial data $\big(X_\red,\{0\}\oplus{\R^r}^*\big)$ relatively to
the smooth family of subspaces $\Delta_\red$;


\end{itemize}

The following facts are immediate:
\begin{enumerate}

\item\label{itm:1} $\Kcal_\Dcal\subset\Kcal_\Delta$ and $\Scal_\Delta\subset\Scal_\Dcal$;

\item\label{itm:2} $\lambda\big(\Scal_{\Delta_\red}\big)=\Scal_\Delta$;

\item\label{itm:3} $\Delta_\red$ is a nondegenerate family of subspaces for
$X_\red$.

\end{enumerate}

We prove the following preparatory lemma:
\begin{lem}\label{thm:conjectura}
An absolutely continuous map $f:[a,b]\to\R^r$ is a solution of $X_\red$ along $\Delta_\red$ iff $v=\sum_{i=1}^rf_iY_i$ is
a solution of
$X$ along $\Delta$. In particular, $\lambda\big(\Kcal_{\Delta_\red}\big)=\Kcal_\Delta\cap\Scal_\Dcal$.
\end{lem}
\begin{proof}
The map $f$ is a solution of $X_\red$ along $\Delta_\red$ iff
$(\mathfrak Bf')_i$ is absolutely continuous for $i=1,\ldots,k$
and
\[\big[\big(\mathfrak Bf'+\mathcal Af\big)_i\big]'=\big(\mathcal Cf+\mathcal A^*f'\big)_i,\quad i=1,\ldots,k.\]
Using \eqref{eq:alphafv} one easily checks that this is also the
condition for $v$ to be a solution of $X$ along $\Delta$.
\end{proof}


We can now prove the main theorem of the section.
\begin{teo}[generalized index theorem]\label{thm:INDEXTHMNOVO}
Consider a triple $(X,\ell_0,\Dcal)$ where $(X,\ell_0)$ is a
symplectic differential system with initial data in $\R^n$ over
the interval $[a,b]$ such that $\ell_0$ defines a nondegenerate
initial condition for $X$ and $\Dcal$ is a smooth family of
subspaces in $\R^n$ over $[a,b]$ whose index with respect to $X$
equals the index of $X$. Denote by $X_\red$ the reduced symplectic
system corresponding to $X$ and $\Dcal$ (see
Remark~\ref{thm:mudaframe}). If $t=b$ is neither a focal instant
for $(X,\ell_0)$ nor a conjugate instant for $X_\red$ then:
\begin{equation}\label{eq:INDEXTHM}
n_-\big(I\vert_{\Kcal_\Dcal}\big)=\iMaslov(X,\ell_0)-\iMaslov(X_\red)+n_-\big(B(a)^{-1}\vert_P\big).
\end{equation}
\end{teo}
\begin{proof}
Since the index of $\Dcal$ with respect to $X$ equals the index of
$X$ then $\Dcal$ is nondegenerate for $X$ and there exists a
maximal negative family $\Delta$ for $X$ with
$\Delta\subset\Dcal$. We may choose a frame $(Y_i)_{i=1}^r$ of
$\Dcal$ such that $(Y_i)_{i=1}^k$ is a frame for $\Delta$, so that
we are in the setup specified at the beginning of the section.
Obviously $\Delta_\red$ is a maximal negative family for $X_\red$,
so that we may apply Theorem~\ref{thm:INDEXvelho} to the triple
$(X,\ell_0,\Delta)$ and to the triple
$(X_\red,\{0\}\oplus{\R^r}^*,\Delta_\red)$ obtaining:
\begin{equation}\label{eq:naoprecisaseparar}
\begin{split}
\iMaslov(X,\ell_0)&\,=n_-\big(I\vert_{\Kcal_\Delta}\big)-n_+\big(I\vert_{\Scal_\Delta}\big)-n_-\big(B(a)^{-1}\vert_P\big),\\
\iMaslov(X_\red)&\,=n_-\big(I_\red\vert_{\Kcal_{\Delta_\red}}\big)-n_+\big(I_\red\vert_{\Scal_{\Delta_\red}}\big).
\end{split}
\end{equation}
From Lemma~\ref{thm:identifyindexformred} we get:
\begin{align}
\label{eq:seila1}n_-\big(I_\red\vert_{\Kcal_{\Delta_\red}}\big)&\,=n_-\big(I\vert_{\lambda(\Kcal_{\Delta_\red})}\big),\\
\label{eq:seila2}n_+\big(I_\red\vert_{\Scal_{\Delta_\red}}\big)&\,=n_+\big(I\vert_{\lambda(\Scal_{\Delta_\red})}\big);
\end{align}
Using item~\eqref{itm:2} on page~\pageref{itm:2}, we have:
\begin{equation}\label{eq:seila3}
n_+\big(I\vert_{\lambda(\Scal_{\Delta_\red})}\big)=n_+\big(I\vert_{\Scal_\Delta}\big).
\end{equation}
Lemmas~\ref{thm:KdSdorthog}, \ref{thm:HKmaisS} and
item~\eqref{itm:1} on page~\pageref{itm:1} imply that we have an
$I$-orthogonal direct sum decomposition
$\Kcal_\Delta=\Kcal_\Dcal\oplus\big(\Scal_\Dcal\cap\Kcal_\Delta\big)$,
where
$\Scal_\Dcal\cap\Kcal_\Delta=\lambda\big(\Kcal_{\Delta_\red}\big)$
by Lemma~\ref{thm:conjectura}. Hence:
\begin{equation}\label{eq:seila4}
n_-\big(I\vert_{\Kcal_\Delta}\big)=n_-\big(I\vert_{\Kcal_\Dcal}\big)+n_-\big(I\vert_{\lambda(\Kcal_{\Delta_\red})}\big).
\end{equation}
The conclusion now follows from equalities \eqref{eq:naoprecisaseparar}---\eqref{eq:seila4}.
\end{proof}

\end{section}

\begin{section}{Geodesics in semi-Riemannian Manifolds}\label{sec:semiRiem}

Let $(M,\gfrak)$ be an $n$-dimensional semi-Riemannian manifold,
where $\gfrak$ is a (nondegenerate) metric tensor of index $k$. We
denote by $\nabla$ the Levi-Civita connection of $\gfrak$ and by
$\Rcal(X,Y)=\nabla_X\nabla_Y-\nabla_Y\nabla_X-\nabla_{[X,Y]}$ the
curvature tensor of $\nabla$. Given a geodesic $\gamma:[a,b]\to M$
then the Jacobi equation along $\gamma$
\begin{equation}\label{eq:Jacobi}
\vfrak''=\Rcal(\gamma',\vfrak)\gamma'
\end{equation}
produces a Morse--Sturm system of the form \eqref{eq:MSJacobi} in
$\R^n$ by means of a parallel trivialization of the tangent bundle
of $M$ along $\gamma$ (the bilinear form $g\in\Bils(\R^n)$
corresponds to the metric tensor $\mathfrak g$ and the
endomorphism $R(t)\in\Lin(\R^n,\R^n)$ corresponds to
$\Rcal(\gamma'(t),\cdot)\gamma'(t)$). In \eqref{eq:Jacobi} the
prime denotes covariant derivative along $\gamma$; this notation
will be used whenever a curve $\gamma$ is fixed by the context.

Different parallel trivializations of $TM$ along $\gamma$ produce
Morse--Sturm systems that are
isomorphic as symplectic differential systems\footnote{%
More generally, non parallel trivializations of $TM$ along
$\gamma$ yield  symplectic differential systems from the Jacobi
equation along $\gamma$. All symplectic differential systems
obtained in this way are isomorphic (see \cite[Section
3]{london}).}.

Let $\Pcal\subset M$ be a smooth submanifold with
$\gamma(a)\in\Pcal$ and $\gamma'(a)\in T_{\gamma(a)}\Pcal^\perp$;
a {\em $\Pcal$-Jacobi field\/} along $\gamma$ is a Jacobi field
$\vfrak$ satisfying the initial conditions:
\begin{equation}\label{eq:PJacobi}
\vfrak(a)\in
T_{\gamma(a)}\Pcal,\quad\gfrak(\vfrak'(a),\cdot)\vert_{T_{\gamma(a)}\Pcal}+\II_{\gamma'(a)}(\vfrak(a),\cdot)=0\in
T_{\gamma(a)}\Pcal^*,
\end{equation}
where $\II_{\gamma'(a)}\in\Bils(T_{\gamma(a)}\Pcal)$ denotes the
second fundamental form of $\Pcal$ in the normal direction
$\gamma'(a)$. If $P\subset\R^n$, $S\in\Bils(P)$ correspond to
$T_{\gamma(a)}\Pcal$ and $\II_{\gamma'(a)}$ by means of the chosen
parallel trivialization of $TM$ along $\gamma$ then $\Pcal$-Jacobi
fields correspond to $(X,\ell_0)$-solutions, where $X$ is the
Morse--Sturm system \eqref{eq:MSJacobi} and
$\ell_0\subset\R^n\oplus{\R^n}^*$ is the Lagrangian
\eqref{eq:ell0PS}. Also, the index form of the pair $(X,\ell_0)$
corresponds to the second variation of the {\em geodesic action
functional\/}
\begin{equation}\label{eq:actionfunc}
E(z)=\frac12\int_a^b\gfrak(z',z')\,\dd t
\end{equation}
at the critical point $\gamma$. The domain of $E$ is the Hilbert manifold $\Omega_{\Pcal q}(M)$ consisting of $H^1$ curves
$z:[a,b]\to M$ with $z(a)\in\Pcal$, $z(b)=q$, where $q=\gamma(b)$. Recall that the critical points of $E$ in $\Omega_{\Pcal
q}(M)$ are the geodesics starting orthogonally at $\Pcal$ and ending at $q$.

The Lagrangian $\ell_0$ defines a nondegenerate initial condition
for $X$ iff the submanifold $\Pcal$ is nondegenerate at
$\gamma(a)$, i.e., if $\gfrak$ is nondegenerate on
$T_{\gamma(a)}\Pcal$. Focal instants for $(X,\ell_0)$ correspond
to $\Pcal$-focal points along $\gamma$. The case where the initial
submanifold $\Pcal$ is a single point corresponds to the case
where $\ell_0=L_0$ (recall \eqref{eq:defL0}); in this case, the
initial condition defined by $\ell_0$ is always nondegenerate.

When $\Pcal$ is nondegenerate at $\gamma(a)$ and $\gamma(b)$ is
not $\Pcal$-focal along $\gamma$ then we can define the {\em
Maslov index\/} $\iMaslov(\gamma,\Pcal)$ of the geodesic $\gamma$
with respect to the initial submanifold $\Pcal$ to be the Maslov
index of the pair $(X,\ell_0)$; by
Proposition~\ref{thm:tudoinvariante}, the Maslov index of $\gamma$
with respect to $\Pcal$ does not depend on the parallel
trivialization used to produce the pair $(X,\ell_0)$. When $\Pcal$
is a single point we call $\iMaslov(\gamma,\Pcal)$ the {\em Maslov
index of $\gamma$\/} and we write simply $\iMaslov(\gamma)$.

If $(\Ycal_i)_{i=1}^r$ are smooth vector fields along $\gamma$
such that $\big(\Ycal_i(t)\big)_{i=1}^r$ is the basis of a
nondegenerate subspace of $T_{\gamma(t)}M$ for all $t\in[a,b]$
then the parallel trivialization along $\gamma$ produce maps
$Y_i:[a,b]\to\R^n$ which form a frame for a nondegenerate smooth
family of subspaces for $X$. The operators $\mathcal A,\mathfrak
B,\mathcal C\in\Lin(\R^r,{\R^r}^*)$ which appear in the
corresponding reduced symplectic system \eqref{eq:sistred} are
given by:
\begin{gather}
\label{eq:novoABCcal}\mathfrak B_{ij}=\gfrak(\Ycal_i,\Ycal_j),\quad\mathcal A_{ij}=\gfrak(\Ycal_j',\Ycal_i),\\
\notag\mathcal C_{ij}=\gfrak(\Ycal_i',\Ycal_j')+\gfrak\big(\Rcal(\gamma',\Ycal_i)\gamma',\Ycal_j\big).
\end{gather}

\begin{defin}\label{thm:reducedMaslov}
Let $\gamma:[a,b]\to M$ be a geodesic and $(\Ycal_i)_{i=1}^r$
smooth vector fields along $\gamma$ such that
$\big(\Ycal_i(t)\big)_{i=1}^r$ is the basis of a nondegenerate
subspace of $T_{\gamma(t)}M$ for all $t\in[a,b]$. Consider the
symplectic differential system $X_\red$ defined in
\eqref{eq:sistred} with $\mathcal A$, $\mathfrak B$, $\mathcal C$
defined in \eqref{eq:novoABCcal}. If $t=b$ is not conjugate for
$X_\red$ then the {\em reduced Maslov index\/} of the geodesic
$\gamma$ (with respect to the fields $\Ycal_i$) is defined by:
\[\iMaslovred(\gamma)=\iMaslov(X_\red).\]
\end{defin}

In this geometrical context, the Index Theorem~\ref{thm:INDEXTHMNOVO} gives a generalized Morse index theorem for
semi-Riemannian geodesics. Observe that the term $n_-\big(B(a)^{-1}\vert_P\big)$ appearing in equality \eqref{eq:INDEXTHM} is
the index of the metric $\gfrak$ in the tangent space $T_{\gamma(a)}\Pcal$ of the initial submanifold.

\subsection{A variational principle for semi-Riemannian geodesics}\label{sub:variational}

We now consider fixed an $n$-dimensional semi-Riemannian manifold
$(M,\gfrak)$ with metric tensor $\gfrak$ of index $k$, a smooth
submanifold $\Pcal\subset M$, a point $q\in M$ and smooth vector
fields $(\Ycal_i)_{i=1}^r$ on $M$ such that
$\big(\Ycal_i(m)\big)_{i=1}^r$ is a basis for a nondegenerate
subspace of $T_mM$ for all $m\in M$. We say that an absolutely
continuous curve $\gamma:[a,b]\to M$ is a {\em geodesic along the
fields $\Ycal_i$\/} if $\gfrak(\gamma',\Ycal_i)$ is absolutely
continuous on $[a,b]$ and
\[\gfrak(\gamma',\Ycal_i)'=\gfrak(\gamma',\Ycal_i'),\]
for $i=1,\ldots,r$. If $\gamma$ is of class $C^2$ then $\gamma$ is a geodesic along the fields $\Ycal_i$ iff $\gamma''$
is orthogonal to the distribution spanned by the $\Ycal_i$; in particular, if
$\gamma$ is a geodesic then
$\gamma$ is a geodesic along the fields
$\Ycal_i$. Moreover, if the vector fields $\Ycal_i$ are Killing\footnote{%
Recall that $\Ycal$ is a Killing vector field iff the bilinear
form $\gfrak(\nabla_\cdot\Ycal,\cdot)$ is skew-symmetric.} then
$\gamma$ is a geodesic along the fields $\Ycal_i$ iff
$\gfrak(\gamma',\Ycal_i)$ is constant for all $i=1,\ldots,r$.

Consider the following subset of the Hilbert manifold $\Omega_{\Pcal q}(M)$:
\[\Ncal_{\Pcal q}(M)=\big\{\gamma\in\Omega_{\Pcal q}(M):\text{$\gamma$ is a geodesic along the fields $\Ycal_i$}\big\}.\]
We are interested in determining conditions that imply that
$\Ncal_{\Pcal q}(M)$ is a Hilbert submanifold of $\Omega_{\Pcal
q}(M)$ and that the critical points of the restriction of the
geodesic action functional $E$ to $\Ncal_{\Pcal q}(M)$ are the
geodesics $\gamma:[a,b]\to M$ starting orthogonally to $\Pcal$ and
ending at $q$. These conditions are given in the following:
\begin{teo}\label{thm:conditions}
Let $\gamma\in\Ncal_{\Pcal q}(M)$ be fixed; consider the following
homogeneous system of linear ODE's in $\R^r\oplus{\R^r}^*$:
\begin{equation}\label{eq:reduzidoestranho}
\left\{\begin{aligned}
&f'=-(\mathfrak B^{-1}\circ\mathcal A)f+\mathfrak B^{-1}\varphi,\\
&(\varphi+\mathcal Ef)'=\big(\mathcal C+\overline{\mathcal
E}-(\mathcal A^*+\mathcal E)\circ\mathfrak B^{-1}\circ\mathcal
A\big)f+\big((\mathcal A^*+\mathcal E)\circ\mathfrak
B^{-1}\big)\varphi,
\end{aligned}\right.
\end{equation}
where $\mathcal A,\mathfrak B,\mathcal C\in\Lin(\R^r,{\R^r}^*)$ are defined in \eqref{eq:novoABCcal} and $\mathcal
E,\overline{\mathcal E}\in\Lin(\R^r,{\R^r}^*)$ are defined by:
\[\mathcal E_{ij}=\gfrak\big(\nabla_{\Ycal_j}\Ycal_i,\gamma'\big),\quad\overline{\mathcal
E}_{ij}=\gfrak\big((\nabla_{\Ycal_j}\Ycal_i)',\gamma'\big).\]
Assume that the system \eqref{eq:reduzidoestranho} does not admit
a non zero solution $(f,\varphi)$ with $f(a)=f(b)=0$. Then:
\begin{enumerate}
\item\label{itm:provar1} $\gamma$ has a neighborhood $\mathfrak V$ in $\Ncal_{\Pcal q}(M)$ which is a Hilbert submanifold of
$\Omega_{\Pcal q}(M)$;

\item\label{itm:provar2} $\gamma$ is a critical point of $E\vert_{\mathfrak V}$ iff $\gamma$ is a geodesic starting
orthogonally at
$\Pcal$;

\item\label{itm:provar3} if $\gamma$ is a critical point of $E\vert_{\mathfrak V}$ then the degeneracy of the second variation
of
$E\vert_{\mathfrak V}$ at $\gamma$ is equal to the multiplicity of $\gamma(b)$ as a $\Pcal$-focal point along $\gamma$.

\end{enumerate}
Assume that $\gamma$ is a critical point of $E\vert_{\mathfrak V}$
such that $\gamma(b)$ is not $\Pcal$-focal along $\gamma$ and
$\Pcal$ is nondegenerate at $\gamma(a)$. If the index of $\gfrak$
restricted to the distribution spanned by $(\Ycal_i)_{i=1}^r$
equals the index of $\gfrak$ then the Morse index of
$E\vert_{\mathfrak V}$ at $\gamma$ is given by:
\begin{equation}\label{eq:TEOREMAO}
n_-\big(\dd^2(E\vert_{\mathfrak
V})(\gamma)\big)=\iMaslov(\gamma)-\iMaslovred(\gamma)+n_-\big(\gfrak\vert_{T_{\gamma(a)}\Pcal}\big).
\end{equation}
\end{teo}

Before getting into the proof of Theorem~\ref{thm:conditions}, which will take up almost entirely the rest of the subsection,
we will make a few remarks about its statement. First we observe that if $\gamma$ is a geodesic then $\overline{\mathcal
E}=\mathcal E'$, so that \eqref{eq:reduzidoestranho} becomes the reduced symplectic system \eqref{eq:sistred}; in particular,
the hypothesis of the theorem is satisfied precisely when $t=b$ is not a conjugate instant for the reduced symplectic system.

Let us look now at the particular case where each $\Ycal_i$ is a
Killing vector field; this obviously implies that $\mathcal
E=-\mathcal A^*$. Another remarkable equality that holds in this
case is $\mathcal C=-\overline{\mathcal E}$; to prove it, recall
that the {\em Hessian\/} of a vector field $\Ycal$ is the
$(2,1)$-tensor field defined by $\Hess(\Ycal)=\nabla\nabla\Ycal$,
i.e.,
$\Hess(\Ycal)(V,W)=\nabla_V\nabla_W\Ycal-\nabla_{\nabla_VW}\Ycal$.
Observe that
$\Hess(\Ycal)(V,W)-\Hess(\Ycal)(W,V)=\Rcal(V,W)\Ycal$; moreover,
if $\Ycal$ is Killing then $\gfrak\big(\Hess(\Ycal)(V,W),Z\big)$
is skew-symmetric in the variables $W$ and $Z$, because
$\gfrak(\nabla_\cdot\Ycal,\cdot)$ is skew-symmetric. Using all
these formulas we compute:
\[\begin{split}
\overline{\mathcal
E}_{ij}&=\gfrak\big((\nabla_{\Ycal_j}\Ycal_i)',\gamma'\big)=\gfrak\big(\Hess(\Ycal_i)(\gamma',\Ycal_j),\gamma'\big)+
\gfrak\big(\nabla_{\Ycal_j'}\Ycal_i,\gamma'\big)\\
&=\gfrak\big(\Rcal(\gamma',\Ycal_j)\Ycal_i,\gamma'\big)+\gfrak\big(\nabla_{\Ycal_j'}\Ycal_i,\gamma'\big)
=-\mathcal C_{ij}.
\end{split}\]

We have proven that if the fields $\Ycal_i$ are Killing then the system \eqref{eq:reduzidoestranho} becomes (recall
\eqref{eq:Aantsim}):
\[\mathfrak Bf'+2\mathcal A_{\mathrm{ant}}f\equiv\text{constant}.\]
Moreover, if the fields $\Ycal_i$ {\em commute\/}, i.e.,
$[\Ycal_i,\Ycal_j]=\nabla_{\Ycal_i}\Ycal_j-\nabla_{\Ycal_j}\Ycal_i=0$
for all $i,j=1,\ldots,r$ then $\mathcal A$ is symmetric and
\eqref{eq:reduzidoestranho} becomes \eqref{eq:veryreduced}.

\begin{rem}\label{thm:KillingCommute}
In the case where the $\Ycal_i$'s are commuting Killing vector
fields then the hypothesis of Theorem~\ref{thm:conditions} is
satisfied iff the symmetric bilinear form $\mathbfB(b)$ (recall
\eqref{eq:Bintegral}) is nondegenerate. Moreover, if $\gamma$ is a
geodesic starting orthogonally at $\Pcal$ then, by\footnote{%
A Killing vector field restricts to a Jacobi field along any
geodesic, so that the fields $\Ycal_i$ indeed correspond to
solutions of $X$.} Remark~\ref{thm:simplyred}, the reduced Maslov
index of $\gamma$ is given by the righthand side of
\eqref{eq:Maslovsimplyred} provided that $\mathfrak B(t)$ is
nondegenerate on $\Img\big(\mathbfB(t)\big)$ for those
$t\in\left]a,b\right]$ such that $\mathbfB(t)$ is degenerate.
\end{rem}

\begin{proof}[Proof of Theorem~\ref{thm:conditions}]
We start by considering the smooth map
\begin{equation}\label{eq:Fcte}
\mathcal F:\Omega_{\Pcal q}(M)\longrightarrow
L^2\big([a,b],{\R^r}^*\big)/\cte
\end{equation}
defined by:
\[\mathcal F(\gamma)(t)_i=\gfrak\big(\gamma'(t),\Ycal_i(\gamma(t))\big)-\int_a^t\gfrak(\gamma',\Ycal_i')\,\dd s\mod{\cte},\]
for all $\gamma\in\Omega_{\Pcal q}(M)$, $t\in[a,b]$ and $i=1,\ldots,r$. In \eqref{eq:Fcte} we have denoted by $\cte$ the
subspace of $L^2\big([a,b],{\R^r}^*\big)$ consisting of constant maps. Obviously:
\begin{equation}\label{eq:NF0}
\Ncal_{\Pcal q}(M)=\mathcal F^{-1}(0).
\end{equation}
The differential of $\mathcal F$ is computed as:
\begin{equation}\label{eq:diffF}
\begin{split}
\dd\mathcal F_\gamma(\vfrak)(t)_i =&\;\mathfrak
g\big(\vfrak'(t),\Ycal_i(\gamma(t))\big)+\gfrak\big(\gamma'(t),\nabla_{\vfrak(t)}\Ycal_i\big)\\
&-\int_a^t
\gfrak(\vfrak',\Ycal_i')+\gfrak\big(\Rcal(\gamma',\vfrak)\gamma',\Ycal_i\big)+\gfrak\big((\nabla_\vfrak\Ycal_i)',\gamma'\big)\,\dd
s\mod{\cte},
\end{split}
\end{equation}
for all $t\in[a,b]$, $i=1,\ldots,r$ and all $\vfrak\in T_\gamma\Omega_{\Pcal q}(M)$. Consider the subspace
$\Scal_\gamma$ of $T_\gamma\Omega_{\Pcal q}(M)$ consisting of vector fields that vanish at the endpoints and that take values
in the span of the fields $\Ycal_i$, i.e.:
\[\Scal_\gamma=\Big\{\sum\nolimits_{i=1}^rf_i\Ycal_i:f_i:[a,b]\xrightarrow{\,H^1\,}\R,\ f_i(a)=f_i(b)=0\Big\}\subset
T_\gamma\Omega_{\Pcal q}(M).\] The central point of the proof is
showing that the restriction of $\dd\mathcal F_\gamma$ to
$\Scal_\gamma$ is an isomorphism; for
$\vfrak=\sum_{i=1}^rf_i\Ycal_i\in\Scal_\gamma$, \eqref{eq:diffF}
can be rewritten as:
\begin{equation}\label{eq:W11L2}
\dd\mathcal F_\gamma(\vfrak)(t)=\mathfrak B(t)f'(t)+\big(\mathcal
A(t)+\mathcal E(t)\big)f(t)-\int_a^t(\mathcal A^*+\mathcal
E)f'+(\mathcal C+\overline{\mathcal E})f\,\dd s\mod{\cte},
\end{equation}
for all $t\in[a,b]$, where $f=(f_i)_{i=1}^r:[a,b]\to\R^r$. The
righthand side of \eqref{eq:W11L2} defines an
$L^2\big([a,b],{\R^r}^*\big)/\cte$-valued Fredholm operator of
index zero in the Hilbert space $H^1_0([a,b],\R^r)$ of $H^1$ maps
$f:[a,b]\to\R^r$ with $f(a)=f(b)=0$. This is proven by an argument
similar to the one used in the proof of Lemma~\ref{thm:HKmaisS},
except that here we also need the compact inclusion
$W^{1,1}\hookrightarrow L^2$. Setting $\varphi=\mathcal
Af+\mathfrak Bf'$ then the righthand side of \eqref{eq:W11L2}
vanishes iff $f$ is a solution of \eqref{eq:reduzidoestranho} with
$f(a)=f(b)=0$; it follows that $\dd\mathcal
F_\gamma\vert_{\Scal_\gamma}$ is injective and therefore an
isomorphism onto $L^2\big([a,b],{\R^r}^*\big)/\cte$.

We can now prove all the assertions made in the statement of the
theorem. Assertion \eqref{itm:provar1} follows from \eqref{eq:NF0}
and from the fact that $\gamma$ is a regular point for $\mathcal
F$. Moreover:
\begin{equation}\label{eq:TNPq}
T_\gamma\mathfrak V=\Ker\big(\dd\mathcal F_\gamma\big).
\end{equation}
Since $\dd\mathcal F_\gamma\vert_{\Scal_\gamma}$ is an
isomorphism, we have:
\[T_\gamma\Omega_{\Pcal q}(M)=T_\gamma\mathfrak V\oplus\Scal_\gamma.\]
Assertion \eqref{itm:provar2} will follow once we establish that $\dd E_\gamma$ vanishes on $\Scal_\gamma$. To see this,
recalling that $\gamma$ is a geodesic along the fields $\Ycal_i$, we compute as follows for
$\vfrak=\sum_{i=1}^rf_i\Ycal_i\in\Scal_\gamma$:
\[\dd E_\gamma(\vfrak)=\int_a^b\mathfrak g(\gamma',\vfrak')\,\dd t=\sum_{i=1}^r\int_a^b\big[f_i\mathfrak
g(\gamma',\Ycal_i)\big]'\,\dd t=0.\]

Assume now that $\gamma$ is a geodesic starting orthogonally at
$\Pcal$. As in the beginning of the section, we choose a parallel
trivialization of $TM$ along $\gamma$ and consider the
Morse--Sturm system with initial data $(X,\ell_0)$ corresponding
to \eqref{eq:PJacobi} and to the Jacobi equation along $\gamma$.
As it was observed, the index form $I\in\Bils(\Hcal)$ of
$(X,\ell_0)$ corresponds to the second variation
$\dd^2E_\gamma\in\Bils\big(T_\gamma\Omega_{\Pcal q}(M)\big)$;
moreover, $\Scal_\gamma$ corresponds to the space $\Scal_\Dcal$ in
\eqref{eq:SD}. Since $\gamma$ is a geodesic, integration by parts
in \eqref{eq:diffF} show that $\mathfrak v\in
T_\gamma\Omega_{\Pcal q}(M)$ is in the kernel of $\dd\mathcal
F_\gamma$ iff $\gfrak(\vfrak',\Ycal_i)$ is absolutely continuous
and
\[\gfrak(\vfrak',\Ycal_i)'=\gfrak(\vfrak',\Ycal_i')+\gfrak\big(\Rcal(\gamma',\vfrak)\gamma',\Ycal_i\big),\]
for all $i=1,\ldots,r$. From \eqref{eq:TNPq} we conclude that the
tangent space $T_\gamma\mathfrak V$ corresponds by the chosen
parallel trivialization of $TM$ along $\gamma$ to the space
$\Kcal_\Dcal$ in \eqref{eq:KD} (here $\Dcal$ is the nondegenerate
family of subspaces for $X$ which has as a frame the maps
$Y_i:[a,b]\to\R^n$ corresponding to the fields $\Ycal_i$). Since
$\gamma$ is a geodesic, the system \eqref{eq:reduzidoestranho}
coincides with the reduced symplectic system $X_\red$, so that
$t=b$ is not conjugate for $X_\red$. The remaining assertions in
the statement of the theorem now follow immediately from
Remark~\ref{thm:KerI} and the generalized Index
Theorem~\ref{thm:INDEXTHMNOVO}.
\end{proof}

\subsection{Geodesics in stationary semi-Riemannian manifolds}\label{sub:stationary}
In this subsection we apply our theory to obtain Morse relations
for geodesics in stationary semi-Riemannian manifolds. For
simplicity, we consider the case of geodesics between two fixed
points. A Ljusternik--Schnirelman theory for this situation was
developed in \cite{GPS-JMAA}; we remark that in this case we will
need only the Index Theorem of \cite{CRAS}
(Theorem~\ref{thm:INDEXvelho}) and not the generalized Index
Theorem~\ref{thm:INDEXTHMNOVO}.

Let $(M,\gfrak)$ be an $n$-dimen\-sional semi-Riemannian manifold
with $\gfrak$ a metric tensor of index $r$. We will call
$(M,\gfrak)$ {\em stationary\/} if it admits Killing vector fields
$(\Ycal_i)_{i=1}^r$ such that $[\Ycal_i,\Ycal_j]=0$ for all
$i,j=1,\ldots,r$ and such that $\big(\Ycal_i(m)\big)_{i=1}^r$ is
the basis of a subspace $\Dcal_m$ of $T_mM$ on which $\gfrak$ is
negative definite for all $m\in M$.

Let $p,q\in M$ be fixed and define $\Ncal_{pq}(M)$ and
$\Omega_{pq}(M)$ as in Subsection~\ref{sub:variational} with
$\Pcal=\{p\}$. Since $\gfrak$ is negative definite on $\Dcal$, the
bilinear form $\mathfrak B(t)\in\Bils(\R^r)$ defined in
\eqref{eq:novoABCcal} is always negative definite and therefore
also $\mathbfB(t)$ is negative definite for all $t\in[a,b]$ (see
\eqref{eq:Bintegral}). It follows from
Remark~\ref{thm:KillingCommute} that the hypothesis of
Theorem~\ref{thm:conditions} is satisfied for every curve
$\gamma\in\Ncal_{pq}(M)$ and that the reduced Maslov index of any
geodesic is zero. Theorem~\ref{thm:conditions} implies the
following facts:
\begin{itemize}
\item $\Ncal_{pq}(M)$ is a Hilbert submanifold of $\Omega_{pq}(M)$;

\item the critical points of $E\vert_{\Ncal_{pq}(M)}$ (see \eqref{eq:actionfunc}) are precisely the
geodesics on $M$ from $p$ to $q$;

\item if $q$ is not conjugate to $p$ then all the critical points of $E\vert_{\Ncal_{pq}(M)}$ are nondegenerate;

\item if $\gamma$ is a nondegenerate critical point of $E\vert_{\Ncal_{pq}(M)}$
then its Morse index is given by:
\[n_-\big(\dd^2(E\vert_{\Ncal_{pq}(M)})(\gamma)\big)=\iMaslov(\gamma).\]

\end{itemize}

\begin{defin}\label{thm:pseudocoercivo}
We say that $E$ is {\em pseudo-coercive\/} on $\Ncal_{pq}(M)$ if given a sequence $(\gamma_n)_{n\ge1}$ in $\Ncal_{pq}(M)$ with
$\sup_{n\ge1}E(\gamma_n)<+\infty$ then $(\gamma_n)_{n\ge1}$ admits a uniformly convergent subsequence.
\end{defin}
Examples and sufficient conditions for $E$ to be pseudo-coercive on $\Ncal_{pq}(M)$ are given in 
\cite[Appendix~B]{GPS-JMAA}.

Let $\Dcal^\perp$ denote the orthogonal complement of $\Dcal$ with respect to $\gfrak$ and let $\gfrak_+$ be the Riemannian
metric in $M$ such that $\Dcal$ and $\Dcal^\perp$ are $\gfrak_+$-orthogonal, $\gfrak_+$ equals $\gfrak$ on $\Dcal^\perp$ and
$\gfrak_+$ equals $-\gfrak$ on $\Dcal$. We define a Riemannian metric on the Hilbert manifold $\Ncal_{pq}(M)$ by:
\[\langle\vfrak,\wfrak\rangle_{H^1}=\int_a^b\gfrak_+(\vfrak',\wfrak')\,\dd t,\quad\vfrak,\wfrak\in
T_\gamma\Ncal_{pq}(M),\ \gamma\in\Ncal_{pq}(M),\]
where the prime denotes the covariant derivative along $\gamma$ with respect to the Levi-Civita connection of $\gfrak_+$.
\begin{prop}\label{thm:PSetcstatio}
If $E$ is pseudo-coercive on $\Ncal_{pq}(M)$ then $E\vert_{\Ncal_{pq}(M)}$ has complete sublevels, it is bounded from below and
it satisfies the Palais--Smale condition. Moreover, if the fields $\Ycal_i$ are complete then $\Ncal_{pq}(M)$ has the same
homotopy type of the loop space of $M$.
\end{prop}
\begin{proof}
See \cite[Proposition 3.3, Theorem 4.1, Proposition 4.3, Proposition 5.2]{GPS-JMAA}.
\end{proof}

\begin{teo}[Morse relations for geodesics in stationary semi-Riemannian manifolds]\label{thm:MORSEstation}
Let $(M,\gfrak)$ be a stationary semi-Riemannian manifold and let $p$ and
$q$ in $M$ be two non conjugate points. For $i\in\N$, set:
\[n_i(p,q)=\text{number of geodesics $\gamma$ in $M$ from $p$ to $q$ with $\iMaslov(\gamma)=i$}.\]
Then, under all the assumptions of Proposition~\ref{thm:PSetcstatio}, we have the following
equality of formal power series in the variable
$\lambda$:
\[\sum_{i=0}^{+\infty}n_i(p,q)\lambda^i=\mathfrak P_\lambda(\Omega^{(0)}(M);\mathbb K)+(1+\lambda)Q(\lambda),\]
where $\mathbb K$ is an arbitrary field, $\Omega^{(0)}(M)$ is the loop space
of $M$, $\mathfrak P_\lambda(\Omega^{(0)}(M);\mathbb K)$ is
its {\em Poincar\'e polynomial\/} with coefficients in $\mathbb K$ and $Q(\lambda)$ is a formal power series in $\lambda$ with
coefficients in $\N\cup\{+\infty\}$.
\end{teo}
\begin{proof}
It follows from Proposition~\ref{thm:PSetcstatio} using standard Morse theory on Hilbert manifolds (see for instance
\cite{oldandnew}).
\end{proof}

\subsection{Geodesics in G\"odel-type spaces}\label{sub:Godel}

In this subsection we will apply our theory to obtain Morse relations for geodesics in semi-Riemannian manifolds of
G\"odel-type; again, we will consider the case of geodesics between two fixed points. A Ljusternik-Schnirelman theory for
this situation (in the Lorentzian case) was developed in \cite{CandelaSanchez}; to obtain the Morse relations we will need the
new index Theorem~\ref{thm:INDEXTHMNOVO} rather than the index theorem of \cite{CRAS}.

Let $(M_0,\gfrak^0)$ be a Riemannian manifold and let $\rho:M_0\to\Bils(\R^r)$ be a smooth map such that $\rho(x)$ is a
nondegenerate symmetric bilinear form of index $k$ in $\R^r$ for all $x\in M_0$. Consider the product $M=M_0\times\R^r$
endowed with the semi-Riemannian metric $\gfrak$ defined by:
\[\gfrak_{(x,u)}\big((\xi_1,\eta_1),(\xi_2,\eta_2)\big)=\gfrak^0_x(\xi_1,\xi_2)+\rho(x)(\eta_1,\eta_2),\]
for all $x\in M_0$, $u\in\R^r$, $\xi_1,\xi_2\in T_xM_0$ and $\eta_1,\eta_2\in\R^r$. In analogy with 
\cite[Definition~1.1]{CandelaSanchez} we will call $(M,\gfrak)$ a {\em semi-Riemannian manifold of
G\"odel type}. In \cite{CandelaSanchez} it is considered the case where $r=2$ and $k=1$.

Consider the commuting Killing vector fields $\Ycal_i=\big(0,\frac\partial{\partial
u_i}\big)$,
$i=1,\ldots,r$ in
$M$. An absolutely continuous curve $\gamma=(\gamma_0,u):[a,b]\to M$ is a geodesic along the fields
$\Ycal_i$ iff
\begin{equation}\label{eq:rhoulinhacte}
\rho\big(\gamma_0(t)\big)u'(t)\equiv\text{constant}\in{\R^r}^*
\end{equation}
for $t\in[a,b]$. Let $p=(p_0,u_0)$, $q=(q_0,u_1)\in M$ be fixed and define $\Ncal_{pq}(M)$ and
$\Omega_{pq}(M)$ as in Subsection~\ref{sub:variational} with
$\Pcal=\{p\}$. For $\gamma=(\gamma_0,u)\in\Ncal_{pq}(M)$, the
bilinear form $\mathfrak B(t)\in\Bils(\R^r)$ corresponding to $\gamma$ defined
in \eqref{eq:novoABCcal} is given by $\mathfrak
B(t)=\rho\big(\gamma_0(t)\big)$; the bilinear form
$\mathbfB(t)\in\Bils({\R^r}^*)$ defined in \eqref{eq:Bintegral} is given by:
\begin{equation}\label{eq:mathbfBrho}
\mathbfB_{\gamma_0}(t)=\int_a^t\rho\big(\gamma_0(s)\big)^{-1}\,\dd s\in\Bils({\R^r}^*),
\end{equation}
where we write $\mathbfB_{\gamma_0}(t)$ rather than $\mathbfB(t)$ to keep the
dependence on $\gamma_0$ explicit. By Remark~\ref{thm:KillingCommute}, the
hypothesis of Theorem~\ref{thm:conditions} is satisfied for $\gamma$ iff
$\mathbfB_{\gamma_0}(b)$ is nondegenerate\footnote{%
In the notation of \cite{CandelaSanchez} (where it is considered the case $r=2$, $k=1$),
the nondegeneracy of \eqref{eq:mathbfBrho} is the condition ``$\vert\mathcal L(x)\vert>0$''.}.

Assume now that for every $\gamma_0\in\Omega_{p_0q_0}(M_0)$ the
bilinear form $\mathbfB_{\gamma_0}(b)$
is nondegenerate. Then a curve $\gamma=(\gamma_0,u)\in\Ncal_{pq}(M)$ is uniquely determined by $\gamma_0$; namely, from
\eqref{eq:rhoulinhacte}, we get:
\begin{equation}\label{eq:ugamma0}
u(t)=u_0+\mathbfB_{\gamma_0}(t)\big(\mathbfB_{\gamma_0}(b)\big)^{-1}(u_1-u_0).
\end{equation}
By Theorem~\ref{thm:conditions}, $\Ncal_{pq}(M)$ is a Hilbert submanifold of
$\Omega_{pq}(M)$; moreover,
we obtain a diffeomorphism $\phi:\Omega_{p_0q_0}(M_0)\to\Ncal_{pq}(M)$ given by $\phi(\gamma_0)=(\gamma_0,u)$ with $u$ defined
in \eqref{eq:ugamma0}. If $E$ denotes the geodesic action functional of $M$ (see \eqref{eq:actionfunc}) then the composite map
$E_0=E\circ\phi:\Omega_{p_0q_0}(M_0)\to\R$ is given by:
\[E_0(\gamma_0)=\frac12\int_a^b\gfrak_0(\gamma_0',\gamma_0')\,\dd t+\frac12\mathbfB_{\gamma_0}(b)^{-1}(u_1-u_0,u_1-u_0),\]
for all $\gamma_0\in\Omega_{p_0q_0}(M_0)$. Theorem~\ref{thm:conditions}
implies that the critical points of $E_0$ are precisely the
curves $\gamma_0\in\Omega_{p_0q_0}(M_0)$ for which $\gamma=\phi(\gamma_0)$ is a geodesic; moreover, $\gamma_0$ is a
nondegenerate critical point of $E_0$ iff $q$ is not conjugate to $p$ along $\gamma$. The index of the second variation of
$E_0$ at a nondegenerate critical point $\gamma_0$ is given by:
\[n_-\big(\dd^2 E_0(\gamma_0)\big)=\iMaslov(\gamma)-\iMaslovred(\gamma).\]
By Remark~\ref{thm:KillingCommute}, the reduced Maslov index
$\iMaslovred(\gamma)$ can be generically computed by formula
\eqref{eq:Maslovsimplyred}.

The Palais-Smale condition and the boundedness from below for the functional $E_0$ are satisfied under certain technical
hypothesis on $\gfrak$. In the result below we will assume that the Hilbert manifold $\Omega_{p_0q_0}(M_0)$ is endowed with
the Riemannian metric:
\[\langle\xi_1,\xi_2\rangle_{H^1}=\int_a^b\gfrak_0(\xi_1',\xi_2')\,\dd t,\quad\xi_1,\xi_2\in
T_{\gamma_0}\Omega_{p_0q_0}(M_0),\ \gamma_0\in\Omega_{p_0q_0}(M_0),\]
where the prime denotes covariant derivative along $\gamma_0$ in the Levi-Civita connection of $(M_0,\gfrak_0)$. Recall that
if $\gfrak_0$ is complete then the metric $\langle\cdot,\cdot\rangle_{H^1}$ is also complete 
(see \cite{Kling}).
\begin{prop}\label{thm:PSetcGodel}
Assume that $(M_0,\gfrak_0)$ is a complete Riemannian manifold, that $\mathbfB_{\gamma_0}(b)$ is nondegenerate for all
$\gamma_0\in\Omega_{p_0q_0}(M_0)$ and that:
\[\sup_{x\in
M_0}\big\Vert\rho(x)^{-1}\big\Vert<+\infty,\quad\sup_{\gamma_0\in\Omega_{p_0q_0}(M_0)}\big\Vert\mathbfB_{\gamma_0}(b)^{-1}\big\Vert<+\infty.\]
Then the functional $E_0:\Omega_{p_0q_0}(M_0)\to\R$ is bounded from below and it satisfies the Palais-Smale condition.
\end{prop}
\begin{proof}
This is proved in \cite[Lemmas 3.5 and 3.7]{CandelaSanchez} in the case $r=2$, $k=1$. The proof of the general case is
analogous.
\end{proof}

The technical hypotheses in the statement of Proposition~\ref{thm:PSetcGodel} are satisfied under suitable boundedness
assumptions on $\rho$ (see \cite[Remark 1.4]{CandelaSanchez} for examples).

\begin{teo}[Morse relations for geodesics in G\"odel-type manifolds]\label{thm:MORSEGODEL}
Let $(M,\gfrak)$ be a semi-Riemannian manifold of G\"odel-type. Let $p=(p_0,u_0)$ and
$q=(q_0,u_1)$ in $M$ be two non conjugate points; for $i\in\N$, set:
\[n_i(p,q)=\text{number of geodesics $\gamma$ in $M$ from $p$ to $q$ with $\iMaslov(\gamma)-\iMaslovred(\gamma)=i$}.\]
Then,
under the assumptions of Proposition~\ref{thm:PSetcGodel}, we have the following equality of formal power series in the
variable
$\lambda$:
\[\sum_{i=0}^{+\infty}n_i(p,q)\lambda^i=\mathfrak P_\lambda(\Omega^{(0)}(M);\mathbb K)+(1+\lambda)Q(\lambda),\]
where $\mathbb K$ is an arbitrary field, $\Omega^{(0)}(M)$ is the loop space
of $M$, $\mathfrak P_\lambda(\Omega^{(0)}(M);\mathbb K)$ is
its {\em Poincar\'e polynomial\/} with coefficients in $\mathbb K$ and $Q(\lambda)$ is a formal power series in $\lambda$ with
coefficients in $\N\cup\{+\infty\}$.
\end{teo}
\begin{proof}
It follows from Proposition~\ref{thm:PSetcGodel} by using standard Morse theory on Hilbert manifolds (see for instance
\cite{oldandnew}) and observing that the loop space of $M$ has the same homotopy type of $\Omega_{p_0q_0}(M_0)$.
\end{proof}

\end{section}



\end{document}